\documentclass[11pt]{article}
\usepackage{amsmath,amsthm,amssymb,amsfonts}
\usepackage{graphicx,color,epsfig,subfigure}
\usepackage{url,verbatim,multirow}
\usepackage[format=plain,  justification=raggedright,singlelinecheck=false, margin=1cm]{caption}

\DeclareMathOperator{\FF}{{\mathcal F}}
\DeclareMathOperator{\GG}{{\mathcal G}}
\DeclareMathOperator{\HH}{{\mathcal H}}

\DeclareMathOperator{\fnd}{{\mathit{f^{(\ell)}(n,d)}}}

\DeclareMathOperator{\cnd}{{\mathit{c^{(\ell)}(n,d)}}}

\DeclareMathOperator{\ex}{ex}
\DeclareMathOperator{\st}{star}
\DeclareMathOperator{\on}{one}
\DeclareMathOperator{\ze}{zero}
\DeclareMathOperator{\cb}{cube}

\newtheorem{thm}{Theorem}
\newtheorem{defn}[thm]{Definition}

\newtheorem{obs}[thm]{Observation}

\newtheorem{clm}[thm]{Claim}

\newenvironment{proofcite}[1]{\noindent{\bf Proof of #1.\,}}{\hfill$\Box$}

\title{On a covering problem in the hypercube }
\author{Lale \"Ozkahya\\
\small \texttt{email: ozkahya@illinoisalumni.org}
\and 
Brendon Stanton\\ 
\small  The Citadel \\
\small  Department of Mathematics  \\
\small  Charleston, SC 29409 USA\\
\small \texttt{email: bstanton@citadel.edu}}

\begin{document}
\maketitle

\begin{abstract}
In this paper, we address a particular variation of the Tur\'an problem 
for the hypercube. 
Alon, Krech and Szab\'o (2007)  asked 
``In an $n$-dimensional hypercube, $Q_n$, and for $\ell<d<n$,   
what is the size of a smallest set, $S$, 
of $Q_\ell$'s so that every 
$Q_d$ contains at least one member of $S$?''  
Likewise, they asked a similar Ramsey type question: 
``What is the largest number of colors that we can use 
to color the copies of $Q_\ell$ in $Q_n$ such that 
each $Q_d$ contains a $Q_\ell$
of each color?''  
We give upper and lower bounds for each of these questions and 
provide constructions of the set $S$ above for some specific cases.  
%for the Tur\'an type problem.
\end{abstract}

\section{Introduction}

For graphs $Q$ and $P$, let
$\ex(Q,P)$ denote the {\it generalized Tur\'an number},
i.e., the maximum number of edges in a $P$-free subgraph of $Q$.
The {\em n-dimensional hypercube}, $Q_n$, is the graph
whose vertex set is $\{0,1\}^n$ and whose edge set is the set of pairs
that differ in exactly one coordinate.
For a graph $G$, we use $n(G)$ and $e(G)$ to denote the number of vertices and
the number of edges of $G$, respectively.

In 1984, Erd{\H o}s~\cite{Erdos} conjectured that
$$ \lim_{n\to \infty}  \frac{\ex(Q_n,C_4)}{e(Q_n)} = \frac{1}{2}.$$
Note that this limit exists, because the function above is non-increasing for $n$ and bounded.
The best upper bound $ \ex(Q_n,C_4)/e(Q_n) \leq 0.62256$ was obtained by
Thomason and Wagner~\cite{Thomason} by
slightly improving the bound 0.62284 given by Chung~\cite{Chung}.
Brass, Harborth and Nienborg~\cite{BHN1995} showed that
the lower bound is
$\frac{1}{2} (1+1/\sqrt{n})$, when $n=4^r$ for integer $r$, and
$\frac{1}{2} (1+0.9/\sqrt{n})$, when $n\geq 9$.

Erd{\H o}s~\cite{Erdos} also asked whether
$o(e(Q_n))$ edges in a subgraph of $Q_n$ would
be sufficient for the existence of a cycle $C_{2k}$ for $k>2$.
The value of $\ex(Q_n,C_6)/e(Q_n)$ is between 1/3 and 0.3941 given by
Conder~\cite{Conder} and Lu~\cite{Lu}, respectively.
On the other hand, nothing is known for the cycle of length 10.
Except $C_{10}$, the question of Erd{\H o}s is answered positively
by showing that $\ex(Q_n,C_{2k}) = o(e(Q_n))$ for $k\ge 4$
in~\cite{Chung}, \cite{Conlon} and \cite{FL2011}.
~\\

A generalization of Erd{\H o}s' conjecture above is
the problem of determining $\ex(Q_n, Q_d)$ for $d\ge 3$.
As for $d=2$, the exact value of $\ex(Q_n,Q_3)$ is still not known.
The best lower bound for $\ex(Q_n,Q_3)/e(Q_n)$ has been
$1 - (5/8)^{0.25} \approx 0.11086$ due to
Graham, Harary, Livingston and Stout~\cite{Graham1993}
until recently Offner~\cite{Offner2011} improved it to 0.1165.
The best upper bound is $\ex(Q_n,Q_3)/e(Q_n)\le 0.25$ due to
Alon, Krech and Szab\'o~\cite{Alon2007}.
They also gave the best bounds for $\ex(Q_n,Q_d)$, $d\ge 4$, as
\begin{equation}\label{ex-qnd}
\Omega(\frac{\log d}{d2^d}) =  1- \frac{\ex(Q_n,Q_d)}{e(Q_n)} \leq
{
\begin{array}{r r}
\frac{4}{(d+1)^2} & \text{if } d \text{ is odd,}\\
\frac{4}{d(d+2)} & \text{if } d \text{ is even.}
\end{array}
}
\end{equation}
These Tur\'an problems are also asked
when vertices are removed instead of edges and
most of these problems are also still open.
In a very recent paper,
Bollob\'as, Leader and Malvenuto~\cite{BLM11} discuss open problems on the vertex-version
and their relation to Tur\'an problems on hypergraphs.
~\\

Here, we present results on a similar dual version of the hypercube Tur\'an problem that is asked by Alon, Krech and Szab\'o in~\cite{Alon2007}.
Let $\HH_n^i$ denote the collection of $Q_i$'s in $Q_n$ for $1\le i\le n-1$.
Call a subset of $\HH_n^\ell$ a {\it $(d,\ell)$-covering set} 
if each member of $\HH_n^d$ contains
some member of this set, i.e., $\HH_n^d$ is {\it covered} by this set. 
A smallest $(d,\ell)$-covering set is called {\it optimal}.
Alon, Krech and Szab\'o~\cite{Alon2007} asked
what the size of the optimal $(d,\ell)$-covering set of $Q_n$ is
for fixed $\ell < d$.
Call this function $\fnd$.
Determining this function when $\ell=1$ is equivalent to the determination
of $\ex(Q_n,Q_d)$, since $\ex(Q_n,Q_d) + f^{(1)}(n,d) = e(Q_n)$ and 
the best bounds for $f^{(1)}(n,d)$ are given in~\cite{Alon2007} 
as~\eqref{ex-qnd}. 
In~\cite{Alon2007}, also the Ramsey version of this problem is asked as follows.
A coloring of $\HH_n^\ell$ is {\em $d,\ell$-polychromatic} if
all colors appear on each copy $Q_d$'s.
%in any $Q\in \HH_n^d$, all colors that are used in this coloring
%appear on copies of $Q_\ell$'s in $Q$.
Let $pc^{(\ell)}(n,d)$ be the largest number of colors for which there exists a
$d,\ell$-polychromatic coloring of $\HH_n^\ell$.
~\\

Let $c^{(\ell)}(n,d)$ be the ratio of $f^{(\ell)}(n,d)$ to the size of $\HH_n^\ell$, i.e.,
\begin{equation}\label{reln-cf}
c^{(\ell)}(n,d) = \frac{f^{(\ell)}(n,d)}{2^{n-\ell}{n\choose \ell}}.
\end{equation}
One can observe that
\begin{equation}\label{reln-cpc}
c^{(\ell)}(n,d) \le \frac{1}{pc^{(\ell)}(n,d)},
\end{equation}
since
any color class used in a $d,\ell$-polychromatic coloring is
a  $(d,\ell)$-covering set of $Q_n$.
Note that the following limits exist, since $c^{(\ell)}(n,d)$ is non-decreasing,
$pc^{(\ell)}(n,d)$ is non-increasing and both are bounded.
$$
c_d^{(\ell)} = \lim_{n\to \infty} c^{(\ell)}(n,d),\quad
p_d^{(\ell)} = \lim_{n\to \infty} pc^{(\ell)}(n,d).
$$
In Section~\ref{section:poly}, we obtain bounds 
on the polychromatic number. 
\begin{thm}\label{polyc}
For integers $n>d>\ell$, let $0<r\le \ell+1$ such that $r = d+1 \pmod{\ell+1}$.
Then
$$
e^{\ell+1} \left(\frac{d+1}{\ell+1}\right)^{\ell+1} \ge {d+1 \choose \ell+1}\ge p_d^{(\ell)} \ge
\left\lceil \frac{d+1}{\ell+1}\right\rceil^r\left\lfloor \frac{d+1}{\ell+1}\right\rfloor^{\ell+1-r}
\approx \left(\frac{d+1}{\ell+1}\right)^{\ell+1}
.$$
\end{thm}
In Section~\ref{section:coverings},
we present the following bounds on $c_d^{(\ell)}$ and $\cnd$.
\begin{thm}\label{covering}
For integers $n>d>\ell$ and $r = d-\ell \pmod{\ell+1}$,
$$
\left(2^{d-\ell}{d\choose \ell}\right)^{-1}
\le c_d^{(\ell)}
\le
\left\lceil \frac{d+1}{\ell+1}\right\rceil^{-r}
\left\lfloor \frac{d+1}{\ell+1}\right\rfloor^{-(\ell+1-r)}.
$$
\end{thm}

% Although the bounds that we provide are very close (NOT TRUE)
The determination of the exact values of $p_d^\ell$ and $c_d^\ell$ 
remains open. 
The lower and upper bounds on $c^{(\ell)}(n,d)$ provided 
in Theorem~\ref{covering} and Theorem~\ref{upbd-f-2}, respectively, 
are a constant factor of each other 
when $d$ and $\ell$ have a bounded difference from $n$.

\begin{thm}\label{upbd-f-2}
Let $n-d$ and $n-\ell$ be fixed finite integers, where $d > \ell$.
%Then, there exists a constant $C = c(n,d,\ell)$ such that
Then, for sufficiently large $n$,
$$
c^{(\ell)}(n,d) \le
%\left\lceil \frac{-(n-d)}{\log{(\frac{(n-d)^{n-d} - (n-d)!}{(n-d)^{n-d}})}}\right\rceil \log{(n-\ell)}
\left\lceil\frac{r\log{(n-\ell)}}{\log{(\frac{r^r}{r^r-r!})}}\right\rceil 
\frac{1+o(1)}{2^{d-\ell}{d\choose l}},
$$
where $r=n-d$.
\end{thm}

Finally, we show an exact result for $\cnd$ when $d=n-1$. 

\begin{thm}\label{lemma:NNMinusOne}
For integers $n-1>\ell$,
$$c^{(\ell)}(n,n-1) =
\frac{\left\lceil \frac{2n}{n-\ell}\right\rceil}{2^{n-\ell}{n\choose \ell}}.
$$
\end{thm}

In our proofs, we make use of the following terminology.
The collection of $i$-subsets of $[n]=\{1,\dots,n\}$, $1\le i\le n$, 
is denoted by ${[n]\choose i}$.  
For an edge $e\in E(Q_n)$, 
$\st(e)$ denotes the coordinate that is different at endpoints of $e$. 
The set of coordinates whose values are 0 (or 1, resp.) at both endpoints of $e$
are denoted by $\ze(e)$ (or $\on(e)$, resp.).
For a subcube $F\subset Q_n$,
$\st(F) := \cup_{e\subseteq E(F)} \st(e)$,
$\on(F) := \cap_{e\subseteq E(F)} \on(e)$ and
$\ze(F) := \cap_{e\subseteq E(F)} \ze(e)$.
Note that $E_1$ {\it covers} $E_2$ for
$E_1\in \HH_n^\ell$ and $E_2\in \HH_n^d$ ($d>\ell$)
if and only if 
$\ze(E_2)\subset \ze(E_1)$ and $\on(E_2)\subset \on(E_1)$.

\begin{defn}\label{wq}
For any $Q\in \HH_n^\ell$ and
$\st(Q)$ with coordinates $s_1<s_2<\ldots < s_\ell$,
we define an $(\ell+1)$-tuple $w(Q)=(w_1,w_2,\dots,w_{\ell+1})$ as\\
- $w_1 = |\{x\in \on(Q): x < s_1\}|$,\\
- $w_j = |\{x\in \on(Q): s_{j-1} < x <s_j\}|$, for $2\le j\le \ell$, \\
- $w_{\ell+1} = |\{x\in \on(Q): x>s_\ell\}|$.\\
\end{defn}

%%%%%%%%%%%%%%%%%%%%%%%%%%%%%%%%%
%%%%%%%%%%%%%%%%%%%%%%%%%%%%%%%%%

\section{Polychromatic Coloring of Subcubes}\label{section:poly}

\begin{proofcite}{Theorem~\ref{polyc}}
{\it The lower bound:}\\
For any $Q\in \HH_n^\ell$ with $w(Q)=(w_1,w_2,\dots,w_{\ell+1})$,
we define the color of each $Q\in \HH_n^\ell$ as the $(\ell+1)$-tuple
$c(Q)=(c_1,\ldots, c_{\ell+1})$ such that
\begin{equation}\label{col-w}
\begin{split}
c_i = w_i \pmod{k}\quad & \text{if $1\le i\le r$ and }\\
c_i = w_i \pmod{k'}\quad & \text{if $r+1\le i\le \ell+1$},  
\end{split}
\end{equation}
where $k = \lceil (d+1)/(\ell+1)\rceil$ and $k'= \lfloor (d+1)/(\ell+1)\rfloor$. 
We show that this coloring is $d,\ell$-polychromatic.

Let $C\in \HH_n^d$, where
$\st(C)$ consists of the coordinates $a_1<a_2<\cdots<a_d$.
We choose a color $(c_1,\ldots, c_{\ell+1})$ arbitrarily and
show that $C$ contains a copy of $Q_\ell$, call it $Q$, with this color.

Since $Q$ must be a subgraph of $C$,
$\ze(C)\subset \ze(Q)$ and $\on(C)\subset \on(Q)$. 
We define $\st(Q) = \{s_1,\ldots,s_\ell\}$ such that
$$
s_i = \begin{cases}
a_{ik} & \text{if $1\le i\le r$,}\\
a_{rk + (i-r)k'} & \text{if $r+1\le i\le \ell$.}
\end{cases}
$$
We include the remaining $d-\ell$ positions of $\st(C)$ 
to $\on(Q)$ or $\ze(Q)$ such that $w(Q)=(w_1,w_2,\dots,w_{\ell+1})$ 
satisfies~\eqref{col-w}. 
This is possible since by the definition of $r$, we have 
$d-\ell = r(k-1) + (\ell+1-r)(k'-1).$ 
%d-\ell = r\left(\left\lceil \frac{d+1}{\ell+1}\right\rceil - 1\right) + 
%(\ell+1-r)\left(\left\lfloor \frac{d+1}{\ell+1}\right\rfloor - 1\right).

{\it The upper bound:}\\
Since $pc^{(\ell)}(n,d)$ is a non-increasing function of $n$,
we provide an upper bound for this function
when $n$ is sufficiently large which is also an upper bound for $p_d^{(\ell)}$.

For a subset $S$ of $[n]$, we define $\cb(S)$ as the subcube $Q$ of $Q_n$
such that $\st(Q) = S$ and $\ze(Q) = [n]\setminus S$. 
Let $\GG$ be a subfamily of $\HH_n^d$ such that
$\GG = \{\cb(S): S\in {[n]\choose d}\}.$
We define a coloring of the members of $\GG$ as follows.

Consider a $d,\ell$-polychromatic coloring of $\HH_n^\ell$ using $p$ colors,
call this coloring $P$.
Fix an arbitrary ordering of the copies of $Q_\ell$'s in $Q_d$.
%By using the coloring $P$ of $\HH_n^\ell$
We define a coloring of $\HH_n^d$ such that
the color of a copy of $Q_d$ is the list of colors of each $Q_\ell$ under $P$ 
in this fixed order. 
By using this coloring on the members of $\GG$, we
obtain a coloring of $\GG$ using $p^{{d\choose \ell}2^{d-\ell}}$ colors.

Now, consider the auxiliary $d$-uniform hypergraph $\GG'$
whose vertex set is the set of coordinates $[n]$ and
whose edge set is defined as the collection of
$\st(E)$'s for each $E$ in $\GG$, i.e.,
$\GG'$ is a complete $d$-uniform hypergraph on the vertex set $[n]$. 
Also we define a coloring of the edges of $\GG'$ by 
using the colors on the corresponding members of $\GG$ as described above. 
Ramsey's theorem on hypergraphs implies that there is a sufficiently large value of $n$
such that
there exists a complete monochromatic subgraph on $d^2+d-1$ vertices
in any edge coloring of $\GG'$ with $p^{{d\choose \ell}2^{d-\ell}}$ colors.
Let $K\subset [n]$ be the vertex set 
of a monochromatic complete subgraph of $\GG'$
on $d^2+d-1$ vertices.
We define $S$ as the collection of $id^{th}$ coordinates in $K$, $1\le i\le d$, 
so that there are at least $d-1$ coordinates between elements of $S$.

\begin{clm}
If $Q$ is a copy of $Q_\ell$ in $\cb(S)$,
then the color of $Q$ under $P$ depends only on $w(Q)$.
\end{clm}
\begin{proof}
Let $E_1$ and $E_2$ be two different copies of $Q_\ell$ in $\cb(S)$ such that
$w(E_1) = w(E_2)$ according to Definition~\ref{wq}.
There exists a subset $S'\subset K$ with $|S'|=d$ such that \\
- $(\on(E_2)\cup \st(E_2)) \subset S'$, i.e., $E_2$ is contained in $\cb(S')$ and\\
- the restriction of $E_2$ on $S'$ gives the same vector as the restriction of $E_1$ on $S$.\\
Clearly, one can find $S'$ that satisfies the first condition. 
It is also possible that $S'$ fulfills the second condition, 
since we can remove or add up to $d-1$ coordinates from $K$ 
between consecutive coordinates of ones and stars in $E_2$ to define $S'$.
This implies that the colors of $E_1$ and $E_2$ are the same under $P$,
since $\cb(S)$ and $\cb(S')$ have the same colors.
\end{proof}

Hence, the number of colors used in
any $d,\ell$-polychromatic coloring of $\HH_n^\ell$
is at most the number of possible vectors $w(Q)$ for any $Q\in \HH_n^\ell$.
The number of possible $(\ell+1)$-tuples $w(Q)$ for any $Q\in \GG$
is given by the number of partitions of at most $d-\ell$ ones into $\ell+1$ parts and
therefore it is at most ${d+1\choose \ell+1}$.
\end{proofcite}

%%%%%%%%%%%%%%%%%%%%%%%%%%%%%%%%%
%%%%%%%%%%%%%%%%%%%%%%%%%%%%%%%%%

\section{The Covering Problem}\label{section:coverings}

\begin{proofcite}{Theorem~\ref{covering}}
Note that a trivial lower bound on $\fnd$ is given by the ratio of $|\HH_n^d|$ 
to the exact number of $Q_d$'s that a single $Q_\ell$ covers in $Q_n$. 
Thus, by~\eqref{reln-cf}, for all $n$,    
\begin{equation}\label{lw-bd}
c^{(\ell)}(n,d)  \ge 
\left\lceil\frac{2^{n-d}{n\choose d}}{{n-\ell \choose n-d}}\right\rceil \cdot 
\frac{1}{2^{n-\ell}{n\choose \ell}}. 
\end{equation}
By using the equality  
${n\choose d}{d \choose d-\ell} = {n\choose \ell}{n-\ell\choose d-\ell}$, 
we are done. 

The upper bound is implied together
by~\eqref{reln-cpc} and Theorem~\ref{polyc}.
\end{proofcite}

We define a {\it $(0,1)$-labelling} of a set as an assignment 
of labels 0 or 1 to its elements.  

\begin{obs}\label{label-cover}
Since any subcube $Q\subset Q_n$ is defined
by $\ze(Q)$ and $\on(Q)$,
a $(d,\ell)$-covering set of $Q_n$ can be defined as
a collection of $(0,1)$-labellings of sets chosen from ${[n]\choose n-\ell}$
such that any $(0,1)$-labelling of sets in ${[n]\choose n-d}$
is contained in at least one of the labelled $(n-\ell)$-sets.
\end{obs}
When providing constructions for the upper bounds
in Theorems~\ref{upbd-f-2} and~\ref{lemma:NNMinusOne},
we provide constructions
for the equivalent covering problem in Observation~\ref{label-cover}.

\begin{proofcite}{Theorem~\ref{upbd-f-2}}
%The lower bounds holds due to Theorem~\ref{covering}. \\

We construct a $(d,\ell)$-covering of $Q_n$  by
providing a construction for the equivalent problem
as stated in Observation~\ref{label-cover}.
In the following, we describe this construction in two steps.
First, we choose the $(n-\ell)$-subsets of $[n]$ to label and 
then, we describe an efficient way to $(0,1)$-label these sets.

{\it Step 1:}
We make use of the following well-known result on the general covering problem.
An {\it $(n,k,t)$-covering} is defined
as a collection of $k$-subsets of $n$ elements
such that every $t$-set is contained in at least one $k$-set.
Let  $C(n,k,t)$ be the minimum number of $k$-sets in an
$(n,k,t)$-covering.
R\"odl proved the following result
by also settling a long-standing conjecture of Erd{\H o}s and Hanani~\cite{EH}.
%\begin{thm}[\cite{Rodl85}]\label{snkt}
For any fixed integers $k$ and $t$ with $2\le t < k < n$,
\begin{equation}\label{snkt}
\lim_{n\to \infty} \frac{C(n,k,t)}{{n\choose t}/{k\choose t}} = 1.
\end{equation}
%\end{thm}
By our assumption, $n-d$ and $n-\ell$ are fixed integers where
$n-d< n-\ell$.
By~\eqref{snkt},
there exists a $(n,n-\ell,n-d)$-covering $\FF$
for sufficiently large $n$ such that
$|\FF|=(1+o(1)){n\choose n-d}/{n-\ell\choose n-d}$.

{\it Step 2:}
We obtain a collection of
$(0,1)$-labellings for each edge $e\in \FF$
so that all $(0,1)$-labellings of $(n-d)$-subsets of $e$ are covered.
%Then, we repeat this procedure for each edge $e$ of $\FF$.
The union of these $(0,1)$-labellings is a covering set. 

%An $r$-uniform hypergraph $\GG$ is called {\it $r$-partite}
%if there are disjoint vertex sets $X_1,\dots,X_r$ such that
%for any edge $E$ of $\GG$, $|E\cap X_i|\le 1$ for $1\le i\le r.$
An {\it r-cut} of an $r$-uniform hypergraph is obtained
by partitioning its vertex set into $r$ parts and taking all edges that meet
every part in exactly one vertex.
An {\it r-cut cover} of a hypergraph is a
collection of $r$-cuts such that each edge is in at least one of the cuts.
An upper bound on the minimum size
of an $r$-cut cover is shown 
by Cioab\u{a}, K\"undgen, Timmons and Vysotsky
in~\cite{Cioaba-Kundgen} 
using a probabilistic proof. 
\begin{thm}[\cite{Cioaba-Kundgen}]\label{mincut}
For every $r$, an $r$-uniform complete hypergraph
on $n$ vertices can be covered with $\lceil c\log{n}\rceil$ 
$r$-cuts if  
$$
c > \frac{-r}{\log{(\frac{r^r - r!}{r^r})}}.
$$
\end{thm}
%\begin{proof}
%Choose a random $r$-cut $R$, where each vertex is placed in one of the $r$ parts with equal probability. Thus, the probability that an edge is not cut by $R$ is $p=(r^r - r!)/r^r$. Assuming that we choose $c\log{n}$ random $r$-cuts, the expected number of edges not cut by any of these cuts is ${n\choose r}p^{c\log{n}}< n^{r+c\log{p}}.$ The bound on $c$ is obtained by letting this number be less than 1.
%\end{proof}

For a fixed edge $e$ of $\FF$,
let $\GG_e$ be the complete $(n-d)$-uniform hypergraph
on the vertex set of $e$.
Let $C=\lceil c\log{(n-\ell)}\rceil$ be the size of a minimum $(n-d)$-cut cover of 
$\GG_e$ as given by Theorem~\ref{mincut}.
We obtain a collection of $(0,1)$-labellings of $e$ 
by labelling each cut in this cover such that 
the vertices in each part are labelled identically with $0$ or $1$. 
Thus, the total number of $(0,1)$-labellings of $e$ is $2^{n-d}C$.
(If some labelling of an edge is used more than once, then we
count this labelling only once.)
Finally, we use similarly labellings for each edge of $\FF$ in the covering set. 
%choose the constant $C$ to be the maximum of its value obtained for each edge.
This yields that
%\begin{equation}\label{cnd-spec}
$$
\cnd
\le \frac{1}{2^{n-\ell} {n\choose \ell}}
(C(1+o(1))2^{n-d}\frac{{n\choose n-d}}{{n-\ell \choose n-d}} =
C(1+o(1))\frac{1}{2^{d-\ell}{d\choose l}},
$$
%\end{equation}
where the last equality is obtained 
%by replacing ${n-\ell \choose n-d}$ with  ${n-\ell \choose d-\ell}$ and
by using the relation
${n\choose d}{d \choose d-\ell} = {n\choose \ell}{n-\ell\choose d-\ell}$.

\end{proofcite}

\begin{proofcite}{Theorem~\ref{lemma:NNMinusOne}}

The lower bound follows from~\eqref{lw-bd}.

For the upper bound, we construct a collection of $(0,1)$-labellings of sets
chosen from ${[n]\choose n-\ell}$,  
where singletons in $[n]$ have both 0 and 1 in some labelling. 
Let $k=\lceil n/ (n-\ell)\rceil$.
We choose a partition $[n]=(P_1,\ldots,P_k)$ such that
$|P_i|=n-\ell$ for $i<k$. 
Let $P\in {[n]\choose n-\ell}$ such that $P_k\subset P$. 
In the covering set, we include two labellings of each of
$P_1,\ldots,P_{k-1},P$, where all labels are the same, either 0 or 1.

%This implies that, if $n-\ell$ divides $n$, then
%$c^{(\ell)}(n,n-1) = \left(2^{n-1-\ell}{n-1\choose \ell}\right)^{-1}.$
\end{proofcite}

~\\


\begin{thebibliography}{20}

\bibitem{Alon2007}
N. Alon, A. Krech and T. Szab\'o.
\newblock Tur\'an's theorem in the hypercube.
\newblock {\em SIAM J. Discrete Math.}, 21: 66--72, 2007.

\bibitem{BLM11}
B. Bollob\'as, I. Leader and C. Malvenuto.
\newblock Daisies and other Tur\'an Problems.
\newblock {\em Combin., Probab. and Comp.}, 20: 743--747, 2011.

\bibitem{BHN1995}
P. Brass , H. Harborth and H. Nienborg.
\newblock On the Maximum number of edges in a $C_4$-free subgraph of $Q_n$.
\newblock {\em J. Graph Theory}, 19: 17--23, 1995.

\bibitem{Chung}
F. Chung.
\newblock Subgraphs of a hypercube containing no small even cycles.
\newblock {\em J. Graph Theory}, 16: 273--286, 1992.

\bibitem{Cioaba-Kundgen}
S. M. Cioab\u{a}, A. K\"undgen, C. M. Timmons and V. V. Vysotsky.
\newblock Covering complete $r$-graphs with spanning complete $r$-partite $r$-graphs.
\newblock {\em Combin., Probab. and Comput.}, 20: 519-527, 2011.

\bibitem{Conder}
M. Conder.
\newblock Hexagon-free subgraphs of hypercubes.
\newblock {\em J. Graph Theory}, 17: 477--479, 1993.

\bibitem{Conlon}
D. Conlon.
\newblock An extremal theorem in the hypercube.
\newblock {\em Electron. J. Combin.}, 17: R111, 2010.

\bibitem{Erdos}
P. Erd\H{o}s.
\newblock On some problems in graph theory combinatorial analysis and combinatorial number theory.
\newblock {\em Graph Theory and Combinatorics}, 1--17, 1984.

\bibitem{EH}
P. Erd{\H o}s and H. Hanani.
\newblock On a limit theorem in combinatorial analysis.
\newblock {\em Publ. Math. Debrecen}, 10: 10--13, 1963.

%\bibitem{EK}
%P. Erd{\H o}s and D. J. Kleitman.
%\newblock On coloring graphs to maximize the proportion of multicolored $k$-edges.
%\newblock {\em J. Combin. Theory}, 5: 164--169, 1968.

%\bibitem{Fredman-Komlos}
%M. L. Fredman and J. Koml\'os.
%\newblock On the size of separating systems and families of perfect hash functions.
%\newblock {\em SIAM J. Alg. Disc. Meth.}, 5: 61--68, 1984.

\bibitem{FL2011}
Z. F\"uredi, L. \"Ozkahya.
\newblock On even-cycle-free subgraphs of the hypercube.
\newblock {\em J. of Combin. Theo., Ser. A}, 118: 1816--1819, 2011.

\bibitem{Lu}
Linyuan Lu.
\newblock Hexagon-free subgraphs in hypercube $Q_n$.
\newblock Private communication.

\bibitem{Offner2011}
D. Offner.
\newblock Some {T}ur\'an type results on the hypercube.
\newblock {\em Discrete Math.}, 102: 2905--2912, 2011.

\bibitem{Offner2008}
D. Offner.
\newblock Polychromatic colorings of the hypercube.
\newblock {\em SIAM J. Discrete Math.}, 22: 450--454, 2008.

\bibitem{Graham1993}
N. Graham, F. Harary, M. Livingston and Q. Stout.
\newblock Subcube fault tolerance in hypercubes.
\newblock {\em Inform. and Comput.}, 102: 280--314, 1993.

\bibitem{Rodl85}
V. R\"odl.
\newblock On a packing and covering problem.
\newblock {\em European J. Combin.}, 6: 69--78, 1985.

\bibitem{Thomason}
A. Thomason and P. Wagner.
\newblock Bounding the size of square-free subgraphs of the hypercube.
\newblock {\em Discrete Math.}, 309: 1730--1735, 2009.

%\bibitem{Wilson1}
%R.M. Wilson.
%\newblock An existence theory for pairwise balanced designs I-II.
%\newblock {\em J. of Combin. Theo., Ser. A}, 13: 220--273, 1972.

%\bibitem{Wilson2}
%R.M. Wilson.
%\newblock An existence theory for pairwise balanced designs III.
%\newblock {\em J. of Combin. Theo., Ser. A}, 18: 71--79, 1975.

\end{thebibliography}
\end{document}